\documentclass{amsart}

\usepackage{amssymb}
\usepackage{url} 
\usepackage{ mathrsfs }

\DeclareMathOperator{\critical}{criticalpt}
\newtheorem{theorem}{Theorem}[subsection]
\newtheorem*{descr}{Description}
\newtheorem{lemma}[theorem]{Lemma}
\newtheorem{prop}[theorem]{Proposition}

\newtheorem{definition}[theorem]{Definition}
\newtheorem{rem}[theorem]{Remark}

\newcommand{\CCC}{{\mathbb C}}

\newcommand{\RR}{{\mathbb R}}

\newcommand{\KK}{{\mathbb K}}

\makeatletter \@addtoreset {equation}{section}

\renewcommand\theequation
  {\ifnum \c@subsection>\z@ \arabic{section}.\arabic{subsection}.\arabic{equation}
  \else \arabic{section}.\arabic{equation} \fi}
\makeatother

\begin{document}


\title[Quantum cohomology of facet symmetric Fano manifolds]{Semisimplicity of the quantum cohomology for smooth Fano toric varieties associated with facet symmetric polytopes.}

\author{Maksim Maydanskiy}
\address{Centre for Mathematical Sciences, Wilberforce Road, Cambridge CB3 0WA
}

\email{M.Maydanskiy@dpmms.cam.ac.uk}

\author{Benjamin P. Mirabelli}

\email{Benjamin.Mirabelli@yale.edu}

\begin{abstract}
The degree zero part of the quantum cohomology algebra of a smooth Fano toric symplectic manifold is determined by the superpotential function, $W$, of its moment polytope.  In particular, this algebra is semisimple, i.e. splits as a product of fields, if and only if all the critical points of $W$ are non-degenerate. 
In this paper we prove that this non-degeneracy holds for all smooth Fano toric varieties with facet-symmetric duals to moment polytopes. \end{abstract}

\maketitle


\section{Introduction}

\subsection{Motivation.}

    In this paper we establish the semisimplicity of the small quantum cohomology algebras of certain smooth Fano toric varieties by studying their associated superpotentials.

    Our main motivation comes from the work of Entov-Polterovich \cite{EnP} and  Ostrover-Tyomkin \cite{OT}, which relates existence of idempotents in the quantum cohomology algebra of a symplectic manifold $(X, \omega)$ to existence of Calabi quasimorphisms and  symplectic quasi-states on this manifold; this in turn has consequences for the group of Hamiltonian symplectomorphisms of $M$ and implies certain restrictions on Lagrangian submanifolds (we refer the reader to Section 1 of \cite{OT} and references therein for further discussion and motivations).  Recently,  \cite{Usher} has generalized these results, showing that semisimplicity of a deformed version of the small quantum cohomlogy is sufficient, and therefore establishing existence of Calabi quasimorphisms and symplectic quasi-states on all compact toric $M$. However, our calculations  of unperturbed quantum cohomology avoid Usher's machinery. Let us also note that different idempotents sometimes produce different quasi-states (see \cite[Corollary F]{OT}, \cite{ElP}; see also \cite[Theorem 1.10]{FOOO} for the deformed version). Recently Borman has used reduction from stable homogeneous quasi-morphisms to construct examples with infinite dimensional space of quasi-states, see \cite{Borm}.  Our family of manifolds with semisimple quantum cohomology (which therefore contains many different idempotents) may lead to other examples of manifolds supporting multiple quasi-states.
    
 \subsection{Quantum cohomology.}  \label{QH}

The version of quantum cohomology we work with is the one appearing in the work of Ostrover-Tyomkin \cite{OT}.  Since there are several different versions of quantum cohomology appearing in the literature, we will now review the definition.

We start with the coefficient ring. We use the Novikov ring $\Lambda^{\scriptscriptstyle
\uparrow}:={\KK^{\scriptscriptstyle \uparrow}}[q,q^{-1}]$ where $\KK^{\scriptscriptstyle \uparrow}$  is the field of generalized Laurent power series,

$$\KK^{\scriptscriptstyle \uparrow} = \Bigl \{
\sum_{\lambda \in \RR}  a_{\lambda} s^{\lambda} \ | \ a_{\lambda}
\in \CCC, \ {\rm and} \ \{ \lambda \, | \, a_\lambda \neq 0 \} \
{\rm is \ discrete \ and \ bounded \ below \ in \ \RR } \Bigr \}.$$

For a $2n$-dimensional symplectic manifold $(X,\omega)$ which is semipositive (for example, Fano), we define the (small) quantum cohomology as $\Lambda^{\scriptscriptstyle \uparrow}$-module to be $QH^*(X, \omega)=H^*(X, \mathbb Q)\otimes_{\mathbb Q} \Lambda^{\scriptscriptstyle \uparrow}$. The quantum cup product is  a deformation of the classical cup product on cohomology, where the structure constants are obtained by counting contributions from certain holomorphic spheres.  The precise definition is as follows ( \cite{McD-S}, Section 11.1): Pick a basis $e_0, \ldots, e_n$ of $H^*(X, \mathbb Q)$ and define a square matrix $q_{\nu \mu}=\int_X e_\nu \cup e_\mu$, and let $q^{\nu \mu}$ be its inverse. Then the quantum cup product is defined by 
\begin{equation}\label{QHProd}
a *b = \sum_{A\in H_2(X)} GW^X_{A,3} (a,b, e_\nu) g^{\nu \mu} e_\mu \otimes s^{\omega (A)} q^{c_1(A)},
\end{equation}
where $GW^X_{A,3} (a,b, c)$ is the 3-pointed Gromov-Witten invariant which, roughly speaking, counts holomorphic spheres passing through Poincare duals of classes $a$, $b$ and $c$.

 \begin{rem} 
 The structure constants of quantum cohomology are obtained by counting contributions from holomorphic spheres as above, and in general, there can be infinitely many spheres contributing.  It is to make sense of the result that one needs to introduce the Novikov ring  $\Lambda^{\scriptscriptstyle
\uparrow}$.  Alternatively, one can take coefficients in a power series ring that tracks all second homology classes that contain contributing holomorphic spheres separately. This is common in algebraic geometry literature (see e.g. \cite{Giv}, \cite{FP}). The smaller coefficient ring $\Lambda^{\scriptscriptstyle
\uparrow}$ is used in part of symplectic literature. For more on various versions of quantum cohomology and their relations, see Section \ref{QHVersions}.
 \end{rem}

 \subsection{Symplectic toric Fano manifolds}
 
 The correspondence between toric varieties and lattice polytopes is well-known (see, for example, \cite{Fulton}, \cite{Cox}, and Sections 2.2, 2.3 and Appendix of \cite{OT}). We recall some basics. 
 
 Symplectically,  a manifold $(X^{2n}, \omega)$ is \emph{toric} if it comes with an effective Hamiltonian action of a $\frac{1}{2} \operatorname{dim}(M)$-dimensional torus $T$. This gives rise to the moment map $\mu: X \mapsto Lie(T)^*$. By a theorem of Atiyah and Guillemin-Sternberg, the image $\mu(X)$ is a polytope, which we will denote  $\check{P}$.  It is defined up to translations in $Lie(T)^*$. Moreover, if $\omega$ is integral, we can pick $\check{P}$ with vertices in the integral lattice of $Lie(T)^*$ and with origin inside the interior of $\check{P}$.

Then  each facet $F$ of  $\check{P}$ is defined by  an inequality $\left<x,n_F \right> \leq a_F$ for $n_F$, the primitive integral inward normal to $F$. The manifold $X$ is \emph{Fano} if and only if all $a_F=1$ for some toric symplectic form $\omega_0$ in the deformation class of $\omega$. A polytope satisfying this condition is called \emph{reflexive}.  The unique special toric form $\omega_0$ as above will be \emph{monotone}, i.e. will satisfy $[\omega_0]=c_1$, where $c_1$ is the first Chern class.

For a polytope $\check{P}$ as above,  the \emph{dual}  polytope $P$  is determined by taking the convex hull of all vectors $\frac{n_F}{a_F}$. This will have vertices on the lattice (in $Lie(T)$) if and only if $\check{P}$ is reflexive.  Clearly, the dual of a reflexive polytope is reflexive. This means that to each symplectic toric Fano manifold we have associated two lattice polytopes - the image of the moment map and its dual. In what follows we will work with the dual of the moment polytope, so that, for example, by ``manifold associated with $P$" we will mean the symplectic toric manifold with moment polytope $\check{P}$.

Finally, we note that given a reflexive polytope $\check{P}$ such that at each vertex~$v$ the set of all normals $n_F$ for  facets $F$ containing~$v$ spans the lattice, it is easy to construct a symplectic toric Fano manifold with moment polytope $\check{P}$ (see \cite{Delz}).

\subsection{Facet symmetric symplectic toric Fano manifolds.}

Reflexive polytopes contain the origin as the only interior lattice point. Indeed, any other interior lattice point  $q$ would lie on a ray to some facet $F$, and we would have $d=\left<q,n_F\right> \in \mathbb{N}$, forcing $a_F>d$ and contradicting $a_F=1$.
By a theorem of Lagarias-Ziegler (\cite{LZ}; see, for example, \cite{Nill}), in each dimension there is, up to a lattice isomorphism, only a finite number of lattice polytopes containing the origin as the only lattice point in the interior. There exists an algorithm for finding all of them, and it produces 4319 isomorphism classes for n = 3, and 473,800,776 for n = 4.
Of course not all of them are smooth, but even the number of smooth ones grows very quickly. The problem of classifying all smooth Fano polytopes has been solved in dimensions less than or equal to 4, but in higher dimensions there is no complete classification (see \cite{Obro} for details).  
 However, adding a mild-looking symmetry assumption changes the situation completely.

\begin{definition}
A polytope $P$ is called \emph{facet-symmetric} if at least one of its facets, $F$, has a reflection with respect to the origin, $-F$, that is also a facet of $P$.
\end{definition} 

 We note that if $P$ is facet-symmetric its dual $\check{P}$ may or may not have the same property, so it is important to keep the two apart.
 
\noindent\textbf{Theorem (\cite{Ewald}):}  Let the vectors $e_i$ be the canonical basis vectors in $\mathbb{R}^n$. A symplectic toric Fano manifolds associated with a facet-symmetric polytope $P$ is a product of some number of copies of
\begin{enumerate}
\item A projective line $\mathbb{P}^1$ with the associated polytope $P=[-1,1] \subset \mathbb{R}$;
 \item The pseudo-del Pezzo varieties $U_k$, where a $U_k$ is of dimension $2k$ and is associated with the polytope $P^k_{pdP}$ that is the convex hull of the vectors $e_1, e_2, \ldots, e_{2k}$, $-e_1, -e_2, \ldots, -e_{2k}$, $e_1+ e_2+\ldots+ e_{2k}$ in $\mathbb{R}^{2k}$;
\item The del Pezzo varieties $V_k$, where a $V_k$ is of dimension $2k$ and is associated with the polytope $P^k_{dP}$ that is the convex hull of the vectors $e_1, e_2, \ldots, e_{2k}$, $-e_1, -e_2, \ldots, -e_{2k}$, $e_1+ e_2+\ldots+ e_{2k}$ and $-e_1- e_2, \ldots-e_{2k}$ in $\mathbb{R}^{2k}$. 

\end{enumerate}

We should note that (despite Ewald's offhand remark to the contrary)  $V_k$ and $U_k$ are \textbf{not} blow-ups of $R=\mathbb{P}^1\times \mathbb{P}^1\times \ldots \times \mathbb{P}^1$ in dimensions $n=2k>2$. Indeed, complex blowing up subtracts $(n-1)$ times the exceptional divisor from the anticanonical class (see \cite{Hartsh}). This corresponds to cutting a corner off of the $n$-cube - the moment polytope of $R$ - and thereby adding a facet $F$ with $a_F=n-1$, and can not produce a Fano toric manifold. This is not such a terrible problem by itself - in the symplectic world changing $a_F$ corresponds to changing the symplectic form of the blow up. But in the case at hand, to obtain $U_k$ from $R$ we chop off a corner at level $a_F=1$, thereby cutting off a number of vertices of the cube, changing its combinatorics and, consequently, the topology of the resulting symplectic manifold. Despite that, $U_k$ and $V_k$ are indeed (smooth) symplectic Fano toric manifolds, as is easy to check from the definition of their associated polytopes.

 This gives a particularly nice family of symplectic toric Fano manifolds.

 \subsection{Semisimplicity of quantum cohomology}
 
 The main result of this paper is

\begin{theorem}\label{Theorem} 
Let $X$ be a facet symmetric smooth symplectic toric Fano manifold. Then for any toric symplectic form $\omega$, $QH^0(X,\omega)$ is semisimple (over the field $\KK^{\scriptscriptstyle \uparrow}$).
\end{theorem}

It is known by work of  Ostrover-Tyomkin that for a Fano toric symplectic manifold  $X$ and a generic toric symplectic form on it, $QH(X,\omega)$ is semisimple. It is also known by work of Entov-Polterovich \cite{EnP} and Ostrover \cite{Ostr} that for a symplectic toric Fano manifold of (real) dimension~4, the $QH(X,\omega)$ is semisimple. (This is an easy consequence of \cite{OT}.) Ostrover-Tyomkin also gives an example (\cite[Section 5]{OT}) of an 8-dimensional Fano toric symplectic manifold for which $QH(X,\omega)$ is \emph{not} semisimple. In general, the question of semisimplicity of $QH(X,\omega)$ for Fano toric symplectic manifolds appears to be a hard combinatorial problem. Theorem  \ref{Theorem} above is a solution of this problem for the special case of facet-symmetric toric Fano manifolds.

\subsection{Versions of quantum cohomology.}\label{QHVersions}

There are different versions of quantum cohomology one can consider. Some of the issues are coefficient rings, gradings and big versus small quantum cohomology. We shall now discuss these variations.

As we explained above, the main point of choosing a coefficient ring for $QH(X,\omega)$  is to make sure that the structure constants in the product \eqref{QHProd} are well-defined. There is a variety of options to choose from. A fairly general construction is described in \cite{McD-S} Definition 11.1.3. (See also \cite{Usher}, Section 7.3.) One starts with a ``3-pointed Gromov-Witten effective cone"  
\begin{equation*}
K^{eff} (X,\omega) \mathrel{\mathop:}= \{ A \in H_2(M) | \exists  A_1, \ldots A_j \in H_2(M) :  A=\sum A_j, GW^X_{A_i,3} \neq 0\}
\end{equation*}

Then a \emph{universal} quantum coefficient ring (over a base field $k$, we use $k= \mathbb{R}$)  is the ring   $$R^0 (X,\omega)  \mathrel{\mathop:}=  \{  \sum_{A \in K^{eff} } \lambda_A e^A , \lambda_A \in R \}$$ of formal sums. These sums can be infinite, but the product is well-defined as a consequence of Gromov compactness. One can restrict consideration to a subring  $$\Gamma (X,\omega) \mathrel{\mathop:}=  \{  \sum  \lambda_A e^A  \in R^0(X, \omega) | \operatorname{sup}_{\lambda(A) \neq 0} c_1(A) < \infty \}.$$ This subring is graded, with $\operatorname{deg} (q^A)= 2 c_1(A)$, the grading appearing in the formula for dimension of moduli space of holomorphic curves used to define $GW^X$.

Finally, a quantum coefficient ring is any ``reduction of coefficients" - namely a ring $\Lambda$ over $k$, with a ring homomorphism  $\phi : R^0(X, \omega) \mapsto \Lambda$, or a graded ring $\Lambda$ with graded ring homomorphism $\phi: \Gamma(X, \omega) \mapsto \Lambda$. We require $\phi$ to preserve the $k$-module structure; finally we want an evaluation map $\iota: \Lambda \mapsto k$, with $\iota(\phi(\sum \lambda_A q^A))=\lambda_0$, and, in the graded case, $\iota(x) =0$ if  $\operatorname{deg} x \neq 0$. 

Then we can define quantum cohomology with coefficients in $\Lambda$. Existence of $\iota$ makes $QH$ a deformation of the regular cup product on $H^*(X)$. In the graded case, the resulting $QH$ is graded (the basis elements $a \in H^*(X)$ are graded by their usual degree in $H^*(X)$). 

Following  \cite{Usher}, we note that all these definitions really depend only on the pair $(X, \mathscr{C})$, where $\mathscr{C}$ is a nonempty connected component of the space of symplectic forms on $X$ (this is because Gromov-Witten invariants are locally constant under variations of the symplectic form). What will depend on the symplectic form, is the reduction of coefficients map $\phi$.

Specifically, given $\omega$, we define the homomorphism from universal quantum coefficient ring to the corresponding Novikov ring $\phi_\omega: \Gamma(X, \omega)\mapsto \Lambda_\omega^{\scriptscriptstyle \uparrow}$ by $\phi_\omega(\sum \lambda_A e^A) = \sum \lambda_A s^{\omega (A)} q^{c_1(A)}$. Defining a quantum cohomology with coefficients in $\Lambda_\omega^{\scriptscriptstyle \uparrow}$ using this $\phi_\omega$, we obtain the definition from Section  \ref{QH}.  The graded piece $QH^0(X, \omega)$ is then a subring and is itself a module over the graded piece of $\Lambda_\omega^{\scriptscriptstyle \uparrow}$ of grading zero, namely the \emph{field} $K^{\scriptscriptstyle \uparrow}$.

Further, in the Fano case, all the sums in equation \eqref{QHProd}  defining structure constants of the quantum cohomology are finite, and one can define a \emph{convergent} version of $QH(X, \mathscr{C})$ with coefficients in $\mathbb{C}$. Namely, a linear map $ev: H_2(X) \mapsto \mathbb{C}$ lets us define for finite sums $EV (\sum \lambda_A e^A)= \sum \lambda_A e^{ev(A)} \in \mathbb{C}$. When $ev$ is of a special form,  $ev(A) = a \omega(A)+b c_1(A)$ for some $\omega \in \mathscr{C}$, $EV$ factors through $\phi_\omega$ (via $s=e^a, q=e^b$). We will denote such maps $EV$ by  $EV_{a,b}$ in this case. In particular, $EV_{a,b}$ can be restricted to $QH^0(X, \omega)$, to give a convergent version of this later ring, which we shall denote $QH^0_{a,b}(X, \omega)$.

There is a $\mathbb{C}^*$ action on the vector space $H^*(X; \mathbb{C})$ underlying all convergent versions $QH^0_{a,b}(X, \omega; \mathbb{C})$, given on $H_l$ by  $z(A)=z^l A$. (This is generated by the so-called Euler vector field, see e.g. \cite[Section 11.5]{McD-S}.)  This action transforms the cup product matrix $g(z)_{ \mu\nu}=z^{2n}  g_{\mu\nu}$,  $g(z)^{\mu\nu}=z^{-2n}  g^{\mu\nu}$, and induces an isomorphism between the convergent quantum cohomologies with $q=e^b$ and $q'=e^b z^2$.  Hence  $QH^0_{a,b}(X, \omega; \mathbb{C})$ only depends on $a$.  Moreover, if $\omega_0$ is monotone Fano, then $ev(A)=  a \omega_0(A)+b c_1(A) = (a+b) c_1(A)$, and all the convergent versions $QH^0_{a,b}(X, \omega; \mathbb{C})$  are isomorphic. In this case $QH^0(X, \omega_0) =QH^0 (X, \omega_0; \mathbb{C}) \otimes K^{\scriptscriptstyle \uparrow}$, and by \cite[Proposition 2.1 ]{EnP}, one is semisimple iff the other is. Hence we can work with the convergent version. We use $QH^0_{0,0} (X, \omega_0; \mathbb{C})$; this convergent version is of course independent of the choice of the form in $\mathscr{C}$. We will call this a \emph{zero-convergent} version of quantum cohomology. It is obtained by setting $q=s=e^0=1$ in the structure constants quantum cohomology from Section \ref{QH}.

 By  \cite[Theorem 4.3]{OT}, if $QH^0(X, \omega_0)$ is semisimple, then $QH^0(X, \omega)$ is for any toric symplectic form $\omega$. Therefore Theorem \ref{Theorem} follows from 

\begin{prop}
Let $X$ be a facet symmetric smooth symplectic toric Fano manifold. Then for monotone symplectic form $\omega_0$, $QH^0(X,\omega_0; \mathbb{C})$ is semisimple (over the field $\mathbb{C}$).
\end{prop}

\subsection{Quantum cohomlogy from superpotentials}

As explained in the previous section, we will study the semisimplicity of zero-convergent quantum cohomology $QH^0(X,\omega_0;\mathbb{C})$. This has a very concrete description as follows. Define the superpotential of a smooth Fano toric symplectic manifold $X$  with the moment polytope $P$ to be the Laurant polynomial given by summing over the facets of the monomials given by the primitive outward normal vectors to those facets, i.e. $W_X=\Sigma_F x^{-n_F}$. Equivalently, $W_X$ is a sum over the vertices of the dual polytope $\check{P}$ of the corresponding monomials.  The Jacobian ideal,  $J$, of $W_X$ in the ring $\mathbb{C}[N]=\mathbb{C}[x_1, x_1^{-1}, \ldots, x_n, x_n^{-1}]$ is generated by all partial \mbox{(log-)derivatives} of $W_X$, and the Jacobian ring is the corresponding quotient. We have \\
 
\begin{descr}If $X$ is a smooth toric Fano symplectic manifold and $W_X$ is as above, then
$QH^0(X) =\mathbb{C}[N]/J$.
\end{descr}

This is the convergent version of \cite[Proposition 3.3]{OT}. 

With this description, the proof of Theorem \ref{Theorem} reduces to showing that the superpotentials associated to facet-symmetric Fano polynomials are non-degenerate. This is what we shall establish in the rest of this paper.

\subsection{Some open problems.}
It is known by \cite[Theorem A]{OT}  that for any Fano symplectic  toric manifold, for generic $\omega$ in $\mathscr{C}$ the quantum cohomology $QH^0(X,\omega)$ is semisimple (over $K^{\scriptscriptstyle \uparrow}$). This means that the convergent versions obtained via $EV_{a,b}$ are semisimple for generic values of $(a,b)$, and in particular that the more general convergent version  is semisimple for generic  $EV$.  It is also known by \cite[Proposition B]{OT},  that there exists $X$ with non-semisimple convergent quantum cohomology for $EV=0$ (eqivalently, non-semisimple quantum cohomology of  monotone $(X, \omega_0)$ over $K^{\scriptscriptstyle \uparrow}$).  However, all known Fano toric manifolds have a field direct summand in quantum cohomology (either in the monotone version over $K^{\scriptscriptstyle \uparrow}$ or, equivalently, in the 0-convergent version over $\mathbb{C}$).  The structure of the monotone/zero-convergent quantum cohomology for general Fano symplectic manifolds is not understood. In particular, it is unknown if there is an $X$ without a field summand in the zero-convergent quantum cohomology, or how to see in general for which $X$ the zero-convergent quantum cohomology is semisimple. One expects that as the dimension goes up, the fraction of toric Fanos with semisimple zero-convergent quantum cohomology approaches 1; a proof of this would likely require a better understanding of reflexive polytopes than is presently available.

An additional set of open questions is that of distinguishing quasimorphisms coming from different idempotents in quantum cohomology. It was shown in \cite[Proposition F]{OT}, that these can be distinct. Another example of this is given in \cite{ElP}. It would be interesting to see if all qusimporphisms obtained from spectral invariants on the facet-symmetric toric $X$ are distinct.

    \subsection{Structure of the paper.}

The rest of this paper is structured as follows. We first determine the superpotential of the product manifold in terms of the superpotential of each of its pieces. Then we investigate the superpotential of each of the component pieces specified by Ewald's theorem. We work directly with polytopes throughout.  We write $W_P$ for the superpotential of manifold $X$ associated with the polytope $P$.  All manifolds are assumed to be toric Fano and smooth.

\subsection*{Acknowledgments}The authors would like to thank John Rickert, Allison Gilmore, David
Jerison, the Massachusetts Institute of Technology, the Center for
Excellence in Education and the Research Science Institute for supporting and funding this research. The first author thanks Yaron Ostrover and Ilya Tyomkin for explanations regarding their work. We thank Alexander Givental and Leonid Polterovich for their comments on earlier versions of this paper. We are also grateful to the referee for useful suggestions and corrections.

\section{Non-degeneracy of the superpotential of
  convex-hull product polytopes.}\label{section convex-hull}

\begin{prop}\label{PROD}
All the critical points of $W_{V_1\times V_2}$ are nondegenerate if and only if all the critical points of $W_{V_1}$ and all the critical points of $W_{V_2}$ are non-degenerate.
\end{prop}

As previously mentioned, taking products of varieties corresponds to taking Cartesian products of their moment polytopes and, hence, to taking the convex-hull product of the corresponding duals.

We investigate the superpotential of $P$ when $P$ is formed by
taking the convex-hull product, $P=\text{Conv}(Q,Q')$.  $\text{Conv}(Q,Q')$ is the convex hull of $(Q\times 0)\cup (0\times Q')$ where
$Q$ is $n$-dimensional and $Q'$ is $m$-dimensional and we use the standard Cartesian product
embedding from $\mathbb{Z}^n$ and $\mathbb{Z}^m$ to $\mathbb{Z}^{n+m}=\mathbb{Z}^n \times \mathbb{Z}^m$.

  If we define $\text{vert}(P)$ to be the set of all vertices of the polytope $P$, then $\text{vert}(P)=(\text{vert}(Q) \times 0) \cup (0 \times \text{vert}(Q')).$  
  
  Hence, the superpotential $W_P$ is the sum of $W_Q$ and $W_{Q'}$, viewed as functions of \emph{different} variables.

Therefore, the set of critical points of $W_P$, $\critical(P)$, is the union of the set of
critical points of $Q\times 0$ and the set of critical points of $0\times Q'$,
$\critical(P)=(\critical(Q) \times 0) \cup (0 \times \critical(Q')).$  If $Q$ is $n$-dimensional and $Q'$ is $m$-dimensional, the Hessian for
$W_P$ is an $(n+m)\times(n+m)$ matrix.  This matrix, $L$, will be of the
form
$$L=\begin{pmatrix}
\text{Hess}(W_Q) & 0\\
0 & \text{Hess}(W_{Q'})
\end{pmatrix}.$$


Thus, $\text{det}(L)\neq 0$ if and only if both $\text{det}(\text{Hess}(W_Q))\neq 0$ and\\ $\text{det}(\text{Hess}(W_{Q'}))\neq 0$.  Therefore, all of the critical points, $(q,q')$, of
$W_P$ are non-degenerate if and only if all the critical points, $q$, of $W_Q$ and all the critical points,
$q'$, of $W_{Q'}$ are non-degenerate.  This proves Proposition~\ref{PROD}.

\subsection{Non-degeneracy of the superpotential of the [-1,1] polytope.}

The one dimensional polytope $[-1,1]$ has superpotential
$W_1=X+\frac{1}{X}$.  The critical points of $W_1$ are the
solutions to the equation $\frac{\partial W_1}{\partial X}=0.$
Therefore, $\pm 1$ are the critical points of $W_1$.  The Hessian of $W_1$ is
$\text{Hess}(W_1)=\left(\frac{\partial ^2W_1}{\partial
  X^2}\right)=\left(\frac{2}{X^3}\right).$ Therefore, $\left.\text{det}(\text{Hess}(W_1))\right\rvert_{\pm 1}=\pm 2\neq 0$.  So,
the superpotential of polytope $[-1,1]$ only has non-degenerate critical points.

\subsection{Non-degeneracy of the superpotential of del Pezzo polytopes.}

We will need the following

\begin{lemma}\label{lem} \ Let 
$\alpha =
\begin{pmatrix}
a & f & f &  \cdots & f\\
f & a & f &  \cdots & f\\
f & f & a &\cdots & f \\
\vdots & \vdots & \vdots& \ddots & \vdots \\
f & f & f &  \cdots & a\\

\end{pmatrix}$ where $\alpha$ is an $L\times L$ matrix.

Let
$\beta =
\begin{pmatrix}
d &  \cdots &d\\
\vdots &  \ddots &\vdots\\
d& \cdots & d\\
\end{pmatrix}$ where $\beta$ is an $(n-L)\times L$ matrix.

Let
$\gamma =
\begin{pmatrix}
d &  \cdots &d\\
\vdots &  \ddots &\vdots\\
d& \cdots & d\\
\end{pmatrix}$ where $\gamma$ is an $L\times (n-L)$ matrix.

Let $\delta=
\begin{pmatrix}
b & h &  h &  \cdots & h\\
h & b & h &  \cdots & h\\
h & h & b &\cdots & h\\
\vdots & \vdots & \vdots& \ddots & \vdots \\
h & h & h &  \cdots & b\\

\end{pmatrix}$ where $\delta$ is an $(n-L)\times (n-L)$ matrix.

The determinant of a $n\times n$ matrix $M$, where $n\geq L\geq 0$, of the form 
$M=\begin{pmatrix}
\alpha & \beta\\
\gamma & \delta
\end{pmatrix}$ is $$\text{det}(M)=(a-f)^{L-1}(b-h)^{n-L-1}((a+f(L-1))(b+(n-L-1)h)-d^2(L(n-L))).$$
\end{lemma}

The proof of this is based on explicitly finding all eigenvectors and eigenvalues of $M$ and is  given in the Appendix.

The superpotential of a del Pezzo polytope $P_{dP}$ of dimension $n$ is
$W_{P_{dP}}=\sum_{j=1}^nX_j+\sum_{j=1}^n\frac{1}{X_j}+\prod_{j=1}^nX_j+\prod_{j=1}^n\frac{1}{X_j}.$
Each critical point of $W_{P_{dP}}$ has coordinates $X_k$ that are the solutions to the set of equations $\frac{\partial W_{P_{dP}}}{\partial X_k}=0$ for all $1\leq
k \leq n$. We compute 
\begin{equation}\label{2}
\frac{\partial W_{P_{dP}}}{\partial
  X_k}=1-\frac{1}{X_k^2}+\prod_{j\neq k}X_j-\frac{1}{X_k^2}\prod_{j\neq
  k}\frac{1}{X_j}=0.
\end{equation}
Multiply \eqref{2} through by $X_k$ to get
\begin{equation}\label{3}
-\left(X_k-\frac{1}{X_k}\right)=\prod_{j=1}^nX_j-\prod_{j=1}^n\frac{1}{X_j}.
\end{equation}
In \eqref{2} and \eqref{3}, $1\leq j\leq n$, where the polytope $P_{dP}$
is $n$-dimensional.
 
 
The product $\displaystyle Z=\prod_{j=1}^nX_j$ is independent of $k$. Therefore,
$\displaystyle\frac{1}{Z}=\prod_{j=1}^n\frac{1}{X_j}.$ Hence, every $X_k$ satisfies the quadratic equation
\begin{equation}\label{quad2}
X_k^2+\left(Z-\frac{1}{Z}\right)X_k-1=0.
\end{equation}

Thus, each coordinate of the critical point is a root of \eqref{quad2}.  Let $A$ and $B$
denote the two roots of \eqref{quad2}. In fact, the roots are $A=-Z$ and $B=\frac{1}{Z}$, but we prefer to keep the $A$ and $B$ as part of our notation for consistency with the later part of the paper.  Therefore, every $X_k$ equals either $A$ or $B$.  
There are $n$ coordinates $X_k$.  Therefore, if $L$ of the $X_k$ coordinates equal $A$ then $(n-L)$ of them equal $B$.  Note $AB=-1$. We have $Z=A^LB^{n-L}=A^L (-1)^{n-L}A^{L-n}=(-1)^L A^{2L-n}$, where we used that $n$ is even. As $Z=-A$, we get 

\begin{equation}\label{A}
(-1)^L A^{2L-n}=-A, 
\end{equation}
\begin{equation}\label{unity}
A^{2L-n-1}=(-1)^{L-1}, 
\end{equation}
so $A$ is an odd root of unity if $L$ is odd and a minus odd root of unity if $L$ is even.

The Hessian of $W_{P_{dP}}$ evaluated at a critical point is of the same form as $M$ in Lemma \ref{lem} because the $X_j$ can only be one of two forms, $L$ of them equal to $B$ and
$n-L$ of them equal to $A$.
From differentiating \eqref{2} and referring back to Lemma \ref{lem}, we see that turning Hess($M$) into the Hess($W_{P_{dP}}$) requires the following
substitutions, 
$$\left. a=\frac{2}{X_1^3}+2X_1^{-(L+2)}X_2^{-(n-L)}\right\rvert_{X_1=A,
  X_2=B},$$
$$\left. b=\frac{2}{X_2^3}+2X_1^{-L}X_2^{-(n-L+2)}\right\rvert_{X_1=A,
  X_2=B},$$
$$\left.
  f=X_1^{L-2}X_2^{n-L}+X_1^{-(L+2)}X_2^{-(n-L)}\right\rvert_{X_1=A,
  X_2=B},$$
$$\left.
  h=X_1^{L}X_2^{n-L-2}+X_1^{-L}X_2^{-(n-L+2)}\right\rvert_{X_1=A,
  X_2=B},$$
$$\left.
  d=X_1^{L-1}X_2^{n-L-1}+X_1^{-(L+1)}X_2^{-(n-L+1)}\right\rvert_{X_1=A,
  X_2=B}.$$

To prove that all critical points of $W_{P_{dP}}$ are non-degenerate we
must prove that $\text{det}(\text{Hess}(M))\neq 0$ and, thus, that $a-f\neq 0$, $b-h\neq 0$, and
$(a+f(L-1))(b+(n-L-1)h)-d^2(L(n-L))\neq 0$, for all critical points of
$W_{P_{dP}}$.

Remembering that $AB=-1$,  that $n$ is even, and that we have equation \eqref{A}, we find
\begin{equation}\label{a}
a=\frac{2}{A^3}+(-1)^L(2A^{n-2L-2})=0,
\end{equation}
\begin{equation}\label{b}
b=-2A^3+(-1)^L(2A^{n-2L+2})=-2A^3-2A,
\end{equation}
\begin{equation}\label{f}
f=(-1)^LA^{2L-n-2}+(-1)^LA^{n-2L-2}=-\frac{1}{A}-\frac{1}{A^3},
\end{equation}
\begin{equation}\label{h}
h=(-1)^LA^{2L-n+2}+(-1)^LA^{n-2L+2}=-A^3-A,
\end{equation}
\begin{equation}\label{d}
d=(-1)^{L-1}A^{2L-n}+(-1)^{L-1}A^{n-2L}=A+\frac{1}{A}.
\end{equation}

So that 
\begin{equation}\label{a-f22}
a-f=\frac{1}{A}+\frac{1}{A^3}
\end{equation}

and 
\begin{equation}\label{b-h22}
b-h=-A^3-A.
\end{equation}

Therefore $a-f$ or $b-h$ cannot equal zero when $A$ is plus or minus an odd root of unity.  Thus, when $a-f$ or $b-h$ is evaluated at a critical point of $W_{P_{dp}}$, $a-f\neq0$ and $b-h\neq0$.

Now we show that $\chi=(a+f(L-1))(b+(n-L-1)h)-d^2(L(n-L))\neq 0$ for any of
the critical points of $W_{P_{dp}}$.  We compute 

$$\chi=(1-L)(\frac{1}{A}+\frac{1}{A^3})(L-n-1)(A+A^3)- L(n-L)(A+\frac{1}{A})^2=(2L-n-1)(A+\frac{1}{A})^2,$$

which is zero at $i$ and $-i$ (remember, $n$ is even), and not at any odd root of unity or its negative. 




Now we have proved that $a-f\neq 0$, $b-h\neq 0$, and
$(a+f(L-1))(b+(n-L-1)h)-d^2(L(n-L))\neq 0$ and, thus, that $\text{det}(M)\neq 0$ for all critical points of
$W_{P_{pdP}}$.  Therefore, all critical points of $W_{P_{pdP}}$ are non-degenerate.


\subsection{Non-degeneracy of the superpotential of pseudo del Pezzo polytopes.}
The superpotential of a pseudo del Pezzo polytope $P_{pdP}$ of
dimension $n$ is $\displaystyle W_{P_{pdP}}=\sum_{j=1}^nX_j+\sum_{j=1}^n\frac{1}{X_j}+\prod_{j=1}^nX_j.$
To find the critical points of $W_{P_{pdP}}$ we must find all
solutions to the set of equations $\frac{\partial W_{P_{pdP}}}{\partial X_k}=0$ for all $1\leq
k \leq n$,
\begin{equation}\label{11}
\frac{\partial W_{P_{pdP}}}{\partial
  X_k}=1-\frac{1}{X_k^2}+\prod_{j\neq k}X_j=0,
\end{equation}
where the polytope $P_{pdP}$ is $n$-dimensional and $1\leq j\leq n.$
Multiply \eqref{11} through by $X_k^2$ to get,
\begin{equation}\label{12}
X_k^2-1+X_k\prod_{j=1}^nX_j=0.
\end{equation}

For any polytope the product $\displaystyle\prod_{j=1}^nX_j$ is independent of $k$. Again using $\displaystyle Z=\prod_{j=1}^nX_j$, every $X_k$ satisfies the quadratic equation
\begin{equation}\label{quad}
X_k^2+ZX_k-1=0.
\end{equation}

Because \eqref{quad} is a second degree equation, $X_k$ has two solutions,
which we again call $A$ and $B$.
From \eqref{quad} we see that $A+B=Z \text{ and } AB=-1.$
Because $\displaystyle Z=\prod_{j=1}^nX_j$ and $X_j$ must equal either $A$ or $B$,
$Z=A^LB^{n-L}$ where $0\leq L\leq n$.  This can be
restated as $Z=(A^LB^L)B^{n-2L}.$  We know that $AB=-1$, so
\begin{equation}\label{C_p}
Z=(-1)^LB^{n-2L}
\end{equation} and 
\begin{equation}\label{AB}
A=-\frac{1}{B}.
\end{equation}
Substitute values for $Z$ and either $A$ or $B$ from \eqref{C_p} and \eqref{AB} into the equation $A+B=Z$ to get
\begin{equation}\label{critpt pdP1}
B-\frac{1}{B}=(-1)^LB^{n-2L}
\end{equation}
and
\begin{equation}\label{critpt pdP2}
A-\frac{1}{A}=(-1)^LA^{2L-n}.
\end{equation}
The Hessian of $W_{P_{pdP}}$ is of the same form as $M$ in Lemma~\ref{lem} because $X_j$ is only one of two forms, with $L$ of them equal to $A$ and
$n-L$ of them equal to $B$.
To make $M$ the Hessian of $W_{P_{pdP}}$ we make the following
substitutions, 
$$\left. a=\frac{2}{X_1^3}\right\rvert_{X_1=A},$$
$$\left. b=\frac{2}{X_2^3}\right\rvert_{X_2=B},$$
$$\left.
  f=X_1^{L-2}X_2^{n-L}\right\rvert_{X_1=A,
  X_2=B},$$
$$\left.
  h=X_1^{L}X_2^{n-L-2}\right\rvert_{X_1=A,
  X_2=B},$$
$$\left.
  d=X_1^{L-1}X_2^{n-L-1}\right\rvert_{X_1=A,
  X_2=B}.$$
To prove that all critical points of $W_{P_{pdP}}$ are non-degenerate we
must prove that $a-f\neq 0$, $b-h\neq 0$, and
$(a+f(L-1))(b+(n-L-1)h)-d^2(L(n-L))\neq 0$ for all critical points of
$W_{P_{pdP}}$. 
Using first \eqref{AB} and then \eqref{critpt pdP2} we find that
\begin{equation}\label{a2}
a=\frac{2}{A^3},
\end{equation}
\begin{equation}\label{b2}
b=-2A^3,
\end{equation}
\begin{equation}\label{f2}
f=(-1)^LA^{2L-n-2}=\left(A-\frac{1}{A}\right)\frac{1}{A^2},
\end{equation}
\begin{equation}\label{h2}
h=(-1)^LA^{2L-n+2}=(A-\frac{1}{A})A^2,
\end{equation}
\begin{equation}\label{d2}
d=(-1)^{L-1}A^{2L-n}=-(A-\frac{1}{A}).
\end{equation}
Now we show that $a-f\neq 0$ and $b-h\neq 0$ for all critical points of $W_{P_{pdp}}$.

From \eqref{a2} and \eqref{f2} we get $a-f=\frac{2}{A^3}-\left(A-\frac{1}{A}\right)\frac{1}{A^2}$. 
From \eqref{b2} and \eqref{h2} we get $b-h=-2A^3-(A-\frac{1}{A})A^2.$  Thus, $a-f=0$ if and only if 
\begin{equation}\label{af}
\frac{3}{A^3}-\frac{1}{A}=0
\end{equation}
 and $b-h=0$ if and only if 
 \begin{equation}\label{bh}
 -3A^3+A=0.
 \end{equation}
 
 By solving \eqref{af} or \eqref{bh} for $A$ we get $A$ is purely imaginary, and substituting into \eqref{critpt pdP2} we see that no possible value solves both \eqref{af} or \eqref{bh} and \eqref{critpt pdP2}.  Therefore, $a-f\neq 0$ and $b-h\neq 0$ for all
critical points of $W_{P_{pdP}}$.

Now we show that $(a+f(L-1))(b+(n-L-1)h)-d^2(L(n-L))\neq 0$ for all critical points of $W_{P_{pdp}}$.
When we make the substitutions from \eqref{a2}-\eqref{d2} we get that
$(a+f(L-1))(b+(n-L-1)h)-d^2(L(n-L))$ is a forth degree
polynomial $U(A)$, with coefficients of $0$ for $A$ and $A^3$,  congruent modulo 2 to $A^4+1$.

 Now by \cite{Ljunggren}, Theorem 3, any polynomial of the same
form as \eqref{critpt pdP2} is irreducible.  Therefore, $U(A)$ must be
divisible by \eqref{critpt pdP2} if it shares any of the same roots. This can happen only when $2L-n$ is 2, 0, or -2.   We computed this polynomial $U(A)$ and the result is 

$$
U(A)=(3-2L-n)A^4+2(2n-5)A^2+(3+2L-3n)$$

but this computation is not necessary, as here we shall use indirect arguments relying only on the above mentioned properties of $U(A)$ to show that none of these cases can occur.

\noindent\textbf{Case 1:} $2L-n=2$, so that $(-1^L)A^2=A-\frac{1}{A}$. Thus, $U(A)= (A^3\pm A^2 \mp 1)(kA+b)$ and, considering the coefficient of $A$, it follows that $k=0$, which is a contradiction.

\noindent\textbf{Case 2:} $2L-n=-2$, so that $(-1^L)A^{-2}=A-\frac{1}{A}$. Thus, $U(A)= (A^3 - A \mp 1)(kA+b)$
and considering the coefficient of $A^3$, it follows that $b=0$, which is a contradiction.

\noindent\textbf{Case 3:} $2L-n=0$, so that $(-1^L)=A-\frac{1}{A}$. Let $U(A)=kA^4+lA^2+m$. Considering $V=V(A)=A^4-3A^2+1=(A^2-A-1)(A^2+A-1)$,  we see that $U(A)- k V(A)= (l+3k)A^2+(m-k)$ is also divisible by $A^2 \pm A -1$, so that $k=m$, $l=-3k$. But we know that $k$ and $m$ are odd and $l$ is even, which implies a contradiction.

Now we have proved that $a-f\neq 0$, $b-h\neq 0$, and
$(a+f(L-1))(b+(n-L-1)h)-d^2(L(n-L))\neq 0$ for all critical points of
$W_{P_{pdP}}$.  Therefore, $\text{det}(W_{P_{pdp}})\neq 0$ when evaluated at any critical point of $W_{P_{pdP}}$ and, thus, all critical points of $W_{P_{pdP}}$ are non-degenerate. 

\subsection{Conclusion}

By \cite{Ewald}, all facet-symmetric smooth Fano polytopes
can be created by taking the convex-hull product of the $[-1,1]$ polytope, del Pezzo polytopes or pseudo del Pezzo
polytopes.  We proved  that the
superpotential of any of these
polytopes only has non-degenerate critical points.  Using the fact that the convex-hull product of any polytopes whose superpotentials only have non-degenerate critical
points creates a polytope whose superpotential only has non-degenerate
critical points we conclude that the superpotential of any facet-symmetric
smooth Fano polytope only has non-degenerate critical points and, thus, the quantum cohomology of its corresponding smooth Fano toric symplectic manifold is semisimple.



\appendix
\renewcommand\thesection{}
\section{}

Here we compute the determinant from Lemma~\ref{lem}.

First, note that the $L \times L$ matrix $\alpha =
\begin{pmatrix}
a & f & f &  \cdots & f\\
f & a & f &  \cdots & f\\
f & f & a &\cdots & f \\
\vdots & \vdots & \vdots& \ddots & \vdots \\
f & f & f &  \cdots & a\\

\end{pmatrix}$ 

\noindent has an eigenvector $v=(1,\ldots, 1)$ with eigenvalue $a+(L-1)f$, and $(L-1)$ eigenvectors orthogonal to $v$ with eigenvalue $(a-f)$. These last $L-1$ eigenvectors, when appended by $n-L$ 0's in the end, survive as eigenvectors of $M$. Similarly the matrix  $\delta$ has an eigenvector $w$ with eigenvalue $b+(n-L-1)h$ and $n-L-1$ eigenvectors orthogonal to $w$ with eigenvalues $b-h$, which when appended with $L$ zeroes in front become eigenvectors of $M$. So we have found $(n-2)$ eigenvectors. What remains is 2 more and, since $M$ is symmetric, they are orthogonal to all the previous ones. Hence, they are of the form $u=kv+lw$ (where we think of $v$ and $w$ as $n$-vectors, by abuse of notation). Plugging in we see that  $Mu=\lambda u$ will hold if and only if $\lambda$ satisfies a quadratic equation with a zeroth degree term $(a+f(L-1))(b+(n-L-1)h)-d^2(L(n-L))$. Hence, the determinant of $M$ is  $(a-f)^{L-1}(b-h)^{n-L-1}((a+f(L-1))(b+(n-L-1)h)-d^2(L(n-L)))$, as advertised.



\end{document}